\documentclass{article}

\usepackage{setspace}
\usepackage{graphicx}
\usepackage[mathscr]{eucal}
\usepackage{enumerate,verbatim}
\usepackage[latin1]{inputenc}
\usepackage{color}
\usepackage{epsf}
\usepackage{graphics}
\usepackage{amsfonts}
\usepackage{amsmath}
\usepackage{epsfig}

\usepackage{authblk}

\begin{document}


\title{Experimental realization of the Furuta pendulum to emulate the Segway motion}


\author[1]{G.Pujol}
\author[1]{L. Acho}
\author[2]{ A. N\'apoles}
\author[3]{V. P\'erez-Garc\'ia} 


\affil[1]{Department of Applied Mathematics III,Universitat Polit\`ecnica de Catalunya-BarcelonaTech (EUETIB), Comte Urgell 187, 08036 Barcelona,
Spain. (email: gisela.pujol@upc.edu, leonardo.acho@upc.edu).}

\affil[2]{Department of Mechanical Engineering, Universitat Polit\`ecnica de Catalunya-BarcelonaTech (EUETIB), Comte Urgell 187, 08036 Barcelona,
Spain. (email: amelia.napoles@upc.edu).}

\affil[3]{Department of Strength of Materials and Structural Engineering, Universitat Polit\`ecnica de Catalunya-BarcelonaTech (EUETIB), Comte Urgell 187, 08036 Barcelona,
Spain. (email: vega.perez@upc.edu).}

\date{}
\maketitle

\doublespacing
\begin{abstract}                          
A Laboratory device is designed to emulate the Segway motion, modifying the Furuta pendulum experiment. To copy the person on the Segway transportation unit to the Furuta pendulum, a second pendulum is added to the main one. Using LMI theory, the control objective consists to manipulate the base of the Furuta pendulum according to the inclination of the added pendulum whereas the coupled pendulums are maintained near to theirs upright positions. According to the inclination of the second pendulum, the base will rotate as long as the second pendulum remains inclined, emulating the Segway behavior. The base rotation will depend on the side of the inclination too.   Experimental results are offered too (video link: http://youtu.be/SHQdW3k3qiE).

\end{abstract}



\section{Introduction}

Under-actuated dynamical systems are those that have more degrees-of-freedom than control inputs. Examples include spacecrafts, underwater autonomous vehicles and mobile robots. Control and stabilization of these systems are challenging tasks and are currently hot topics of research for both engineers and applied mathematicians \cite{HW13,STSK13}. New stabilization strategies are validated and tested on classical benchmark systems such as the `ball on a beam' and `inverted pendulum' systems \cite{AF00,awtar}. An interesting problem comes from introducing some kind of perturbation to these dynamical systems, to study the robustness in front external disturbances \cite{asian12,PA10}. Also, one such systems which has drawn the attention of control researchers is the Mobile Inverted Pendulum (MIP) which is a two-wheeled robot with a central body that carries a payload. The robot has the advantage of having a small footprint in addition to its ability to turn about its central axis. A commercial variant of the MIP is the well known Segway \cite{N2004}. We can define a Segway Robotic Mobile Platform (RMP) as self-balancing personal transport unit which is similar to the classic control inverted pendulum stabilization problem. 
\begin{figure}[ht!]
\begin{center}
 \includegraphics[scale=0.56]{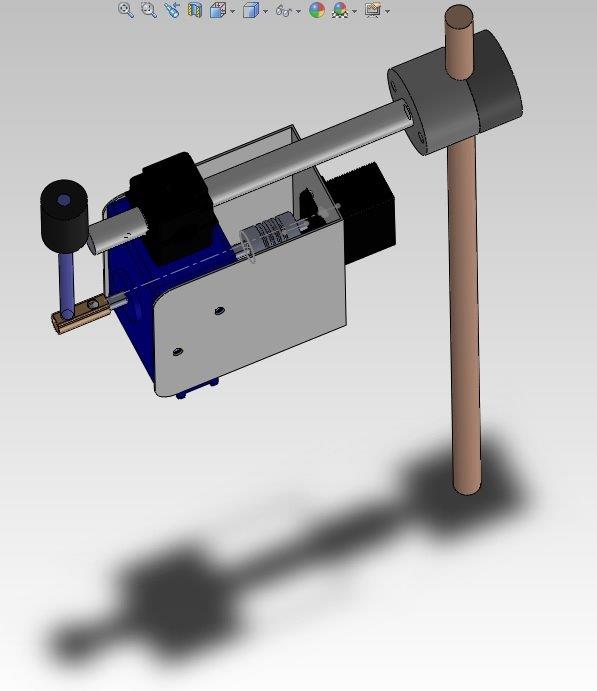}   
\caption{Laboratory device.}  
\label{fig0}                                 
\end{center}                                 
\end{figure}
In this paper, a laboratory device is designed to create a commanded disturbance to the Furuta pendulum (rotational inverted pendulum). This pendulum is modified by adding a second inverted pendulum coupled to the main one (see Figure \ref{fig0}). The induced motion on the second pendulum causes displacement of the center mass of the system, producing a kind of perturbation similar to that presented on MIP transportation unit. Roughly speaking, modifying the position of an inverted pendulum (the person driving the Segway),  the displacement of the base (the Segway) is obtained. In practice, we modify an existent experiment Furuta pendulum presented in Figure \ref{maq-pendu} (provided by ECP-systems) by replacing weight $m_{w}$ by this second inverted pendulum.  Figures \ref{fig11} and \ref{fig1} present this new experimental realization. This external perturbation (second pendulum) is quantized by measuring its angular position via an encoder from  {\it Passport} company. So, we can study directly how this disturbance actuates on our inverted pendulum. 
\begin{figure}[ht!]
\begin{center}
\includegraphics[height=6cm,width=7cm]{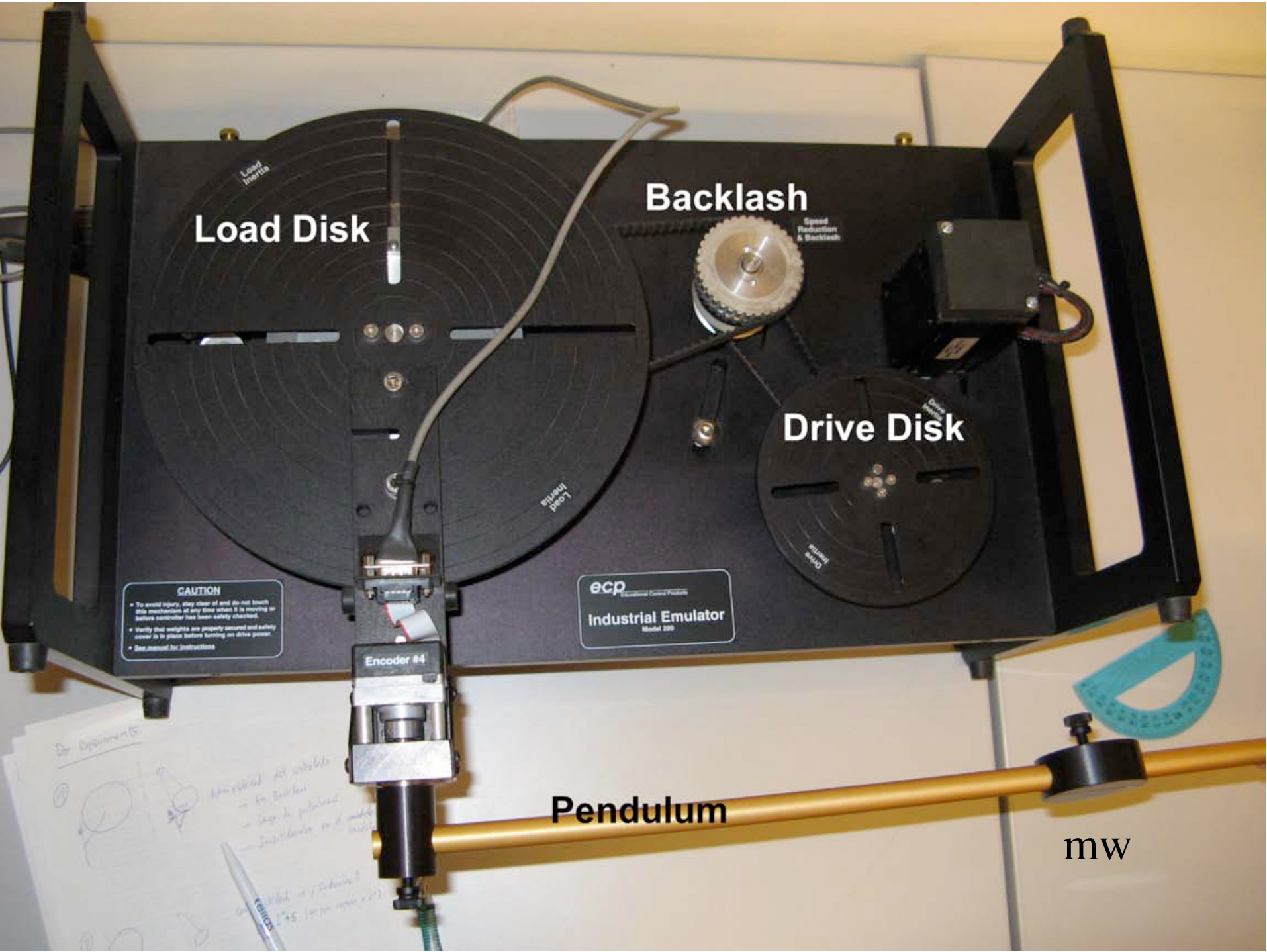}
\caption{Basic pendulum assembly control (Model M220) of ECP systems company \cite{M220}.} \label{maq-pendu}
\end{center}
\end{figure}

Our control objective is to  experiment with this dynamic system in order to prove robust stability of an output feedback controller.  Many nonlinear results exist in the bibliography, incorporating experimental applications \cite{SWS04,TK96,TT03,MS13}, but the purpose of the present paper is to design an easy control algorithm that solves a complex nonlinear problem. For linear systems, $H_\infty$ control theory offers the
possibility of including robustness considerations explicitly in the design and the opportunity to formulate physically meaningful performance objectives that can be expressed as $H_\infty$ design specifications \cite{doyle,KPZ90}, and solved via the  Linear Matrix Inequality (LMI)  techniques \cite{OP06,c1}. It is convenient to point that
firstly the $H_{\infty}$ linear controller clammed in \cite{doyle} was implemented in our Furuta experiment, but experimentally
the closed-loop systems was unstable \cite{PA10}. Based on LMI techniques, \cite{Aziz10} presents a hierarchy intelligent control scheme for a Segway vehicle, but without real experimentation. Also, a fuzzy control is designed and implemented to the same problem in \cite{HWC11}. With simplicity in mind, we propose a linear approach. An output feedback LMI controller is obtained from the unperturbed model (the one without the second inverted pendulum). The stability is proved, and the controller performance is evaluated experimentally. According to experiments, the output feedback LMI controller is robust against this kind of perturbation. In particular,  its possible realization as a Segway navigation control unit center is proved. For control application on the Segway, we use a LMI-$H_\infty$ controller originally designed to stabilize the inverted pendulum. We modify it to fit the navigation control objective. This change consists to do not take into account the displacement variable of the load disk (see Figure \ref{maq-pendu}), that is, the output control ignores the disk  position. Doing this modification on the control law,  the disk position moves freely and we obtain navigation.

This paper presents a laboratory device motivating first the Segway realization. Then, mathematical model as well as the control problem are firstly described in Section three. Experimental test and results are commented in Section four, showing controller robustness. Finally, Section five gives the conclusions.

\section{Laboratory device design}

\subsection{Discussion: segway inspiration}
\begin{figure}[ht!]
\begin{center}
 \includegraphics[scale=0.56]{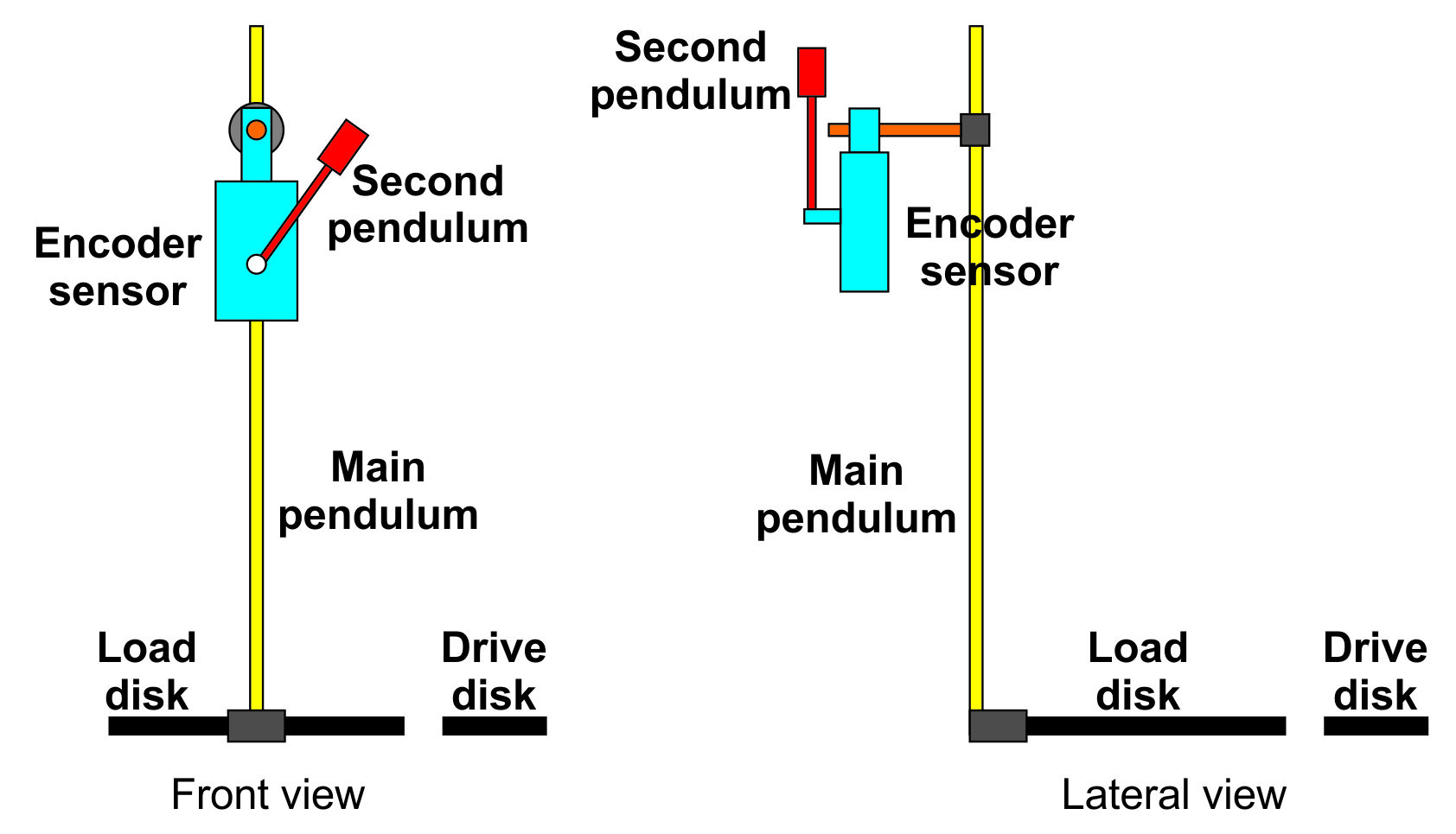}   
\caption{Coupled pendulums experiment: concept design (front and lateral view). On this picture, in blue, the encoder sensor. In yellow, the main pendulum. In red,  the second pendulum coupled to the main one.}  
\label{fig11}                                 
\end{center}                                 
\end{figure}

\begin{figure}[ht!]
\begin{center}
\includegraphics[scale=0.6]{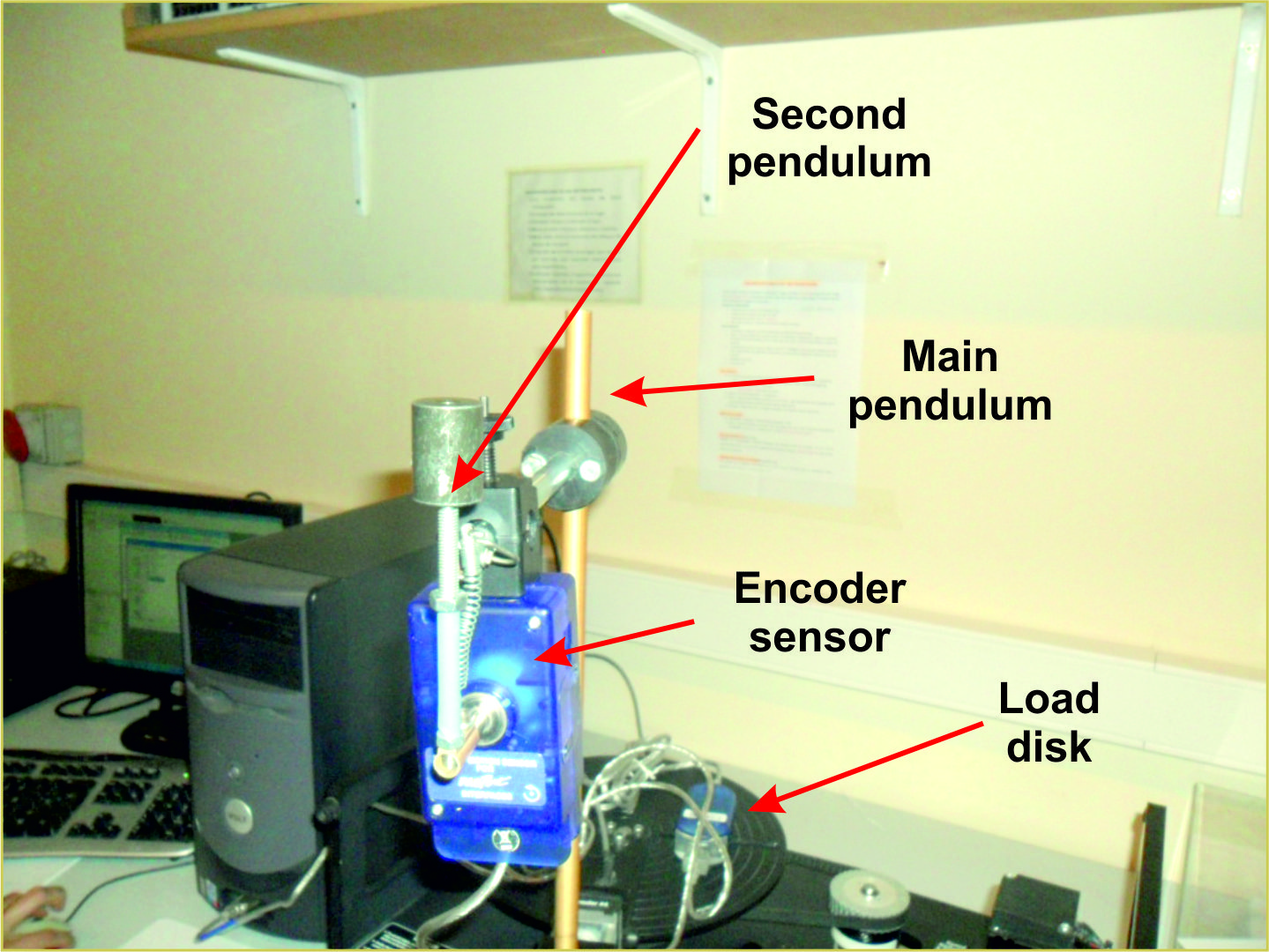} 
\caption{Coupled pendulums experiment: A new realization.}  
\label{fig1}                                 
\end{center}                                 
\end{figure}

As it was said in the introduction, under-actuated dynamical systems are those that have more degrees-of-freedom than control inputs. One such systems which has drawn the attention of control researchers is the MIP. The robot has the advantage of having a small footprint in addition to its ability to turn about its central axis. A commercial variant of the MIP is the well known Segway \cite{N2004}. We can define a Segway Robotic Mobile Platform (RMP) as self-balancing personal transport unit which is similar to the classic control inverted pendulum stabilization problem. The navigation of this unit is due to, basically, the displacement of the center of mass of the body. This line of reasoning motivates us our design showed in Figures \ref{fig11} and \ref{fig1}. With this new device, the position of the center of gravity of the second inverted pendulum is modified by hand, to induce a perturbation which produces displacement of the Furuta pendulum. At this respect, the main pendulum, attached to the rotating base (Figures \ref{maq-pendu}, \ref{fig11}, \ref{fig1} and \ref{fig3}), can emulate the Segway machine, with the difference that the Segway moves linearly, but the Furuta pendulum rotates. The second inverted pendulum can emulate a person (Figure \ref{fig3}).

\begin{figure}[ht!]
\begin{center}
\caption{Basic Segway inspiration scheme.}
\label{fig3}
\end{center}
\end{figure}

\subsection{Description of the laboratory device}
The laboratory device is described in Figure \ref{fig1}. The original experiment is formed by the load and drive disk, with the main pendulum, as Figure \ref{maq-pendu} shows. These parts emulate the Segway. In this experiment, we consider the load and drive disk as an ensemble, and we define it as the load disk. At $y_{m}=32$ cm (Figure \ref{441}), we place the new device laboratory that emulates the person driving the Segway. Its total weight is $48$ gr. The second pendulum is added to change its center of mass. Its weight is $21$ gr. It is possible to modify the weight and structure of the second pendulum changing its end-piece. We have angular position measurements of the disk ($\theta_{1}$) and the main pendulum ($\theta_{2}$). Also, the angular position of the second pendulum is measured off-line from the control-loop via an encoder.

\section{Dynamical system equations and control statement}

\subsection{Dynamical system equations}
The rotating base system driven by a motor shown in Figures \ref{maq-pendu} and \ref{441} is representative of the 
pendulum attached to experimental Model 220 apparatus \cite{M220}, where $\theta_1(t)$ is the magnitude of angular rotation of the disk with respect to $O_1$ and
$\theta_2(t)$ is the magnitude of angular rotation of the pendulum with respect to $O_2$. In the case study described in this paper, the mass set at $P$ (weight $m_{w}$ in Figure \ref{441}) is changed by a second inverted pendulum with total mass $48$gr (Figures  \ref{fig11} and \ref{fig1}). 

\begin{figure}[ht!]
\centering
\caption{Diagram representation of Furuta pendulum with rotating
base \cite{M220}.} 
\label{441}
\end{figure}

The equations of motion for the unperturbed case are obtained from
Lagrange's equations and may be linearized in the case of $\theta_1=0$ and
$\theta_2=0$ \cite{M220}:
\begin{align}
\ddot{\theta}_1 = \frac{1}{p}\left( -\dot{\theta}_1 \bar{J}_z
-m^2l^2_{cg}R_hg\theta_2+\bar{J}_zu\right)\label{ein}\;,\\
\ddot{\theta}_2 = \frac{1}{p}\left( mR_hl_{cg}\dot{\theta}_1 -
mgl_{cg}(\bar{J}_1+J_y)\theta_2-mR_hl_{cg}u\right)\;,\label{zwei}
\end{align}
\noindent where $p=\bar{J}_z(\bar{J}_1+J_y)-(mR_hl_{cg}^2)$ is normalized,
$\bar{J}_1=J_1+mR_h^2$ and $\bar{J}_z=J_z+ml_{cg}\,$; $J_1$
includes equivalent inertia of all elements that move uniformly
with the motor disk; $J_y$ and $J_z$ are the pendulum moments of inertia, relative to its center of mass; $l_{cg}$ is the center of gravity of the combined pendulum road and weight ($m=m_{r}+m_{w}$); $g$ is the
gravity constant; and $u(t)$ is the control effort. For more
details on parameters and modeling, see \cite{M220}.

As equations (\ref{ein}) and (\ref{zwei}) suggest, the state vector is defined by
\begin{equation}
\label{xx}\mathbf{x}^T(t)=[\theta_1(t)\;\dot{\theta}_1(t)\;\theta_2(t)\;\dot{\theta}_2(t)]\;.
\end{equation}
\noindent Our control purpose is to navigate the load disk, induced by the second pendulum. This scenario is defined considering the load disk plus the main pendulum as the Segway, and the second pendulum as the person driving it. So, the feedback control must consider the position $\theta_2(t)$ and velocity $\dot{\theta}_2(t)$ of the main pendulum induced by the second pendulum (by hand in the experiment, but in the future modified using radio control). Also, the velocity of the load $\dot{\theta}_1(t)$ has to be taken into account, because the velocity of the Segway must be controlled.

A linear output feedback strategy is adopted to simplify the approach. In fact, it is convenient to point that firstly the $H_{\infty}$ linear controller clammed in \cite{doyle} was implemented in our Furuta experiment. Notwithstanding, the experimental test indicates that the closed-loop system was unstable. Also, when the state-feedback control was implemented (full information case), the disk returned to its equilibrium position, and we did not obtain navigation (or Segway displacement). To obtain navigation, we have to ignore the load disk position $\theta_{1}(t)$. That is the reason to consider an output feedback control defined by $
u(t)=\mathbf{K}\mathbf{y}(t)
$, with $\mathbf{y}^T(t)=[\dot{\theta}_1\;\theta_2\;\dot{\theta}_2]$. In matrix notation, we have:
\begin{equation}
\label{out}
\mathbf{y}(t)=\underbrace{\left[\begin{array}{cccc}
0&1&0&0\\[2ex]
0&0&1&0\\[2ex]
0&0&0&1\\[2ex]
\end{array}\right]}_{\mathbf{C}_2}\mathbf{x}(t)
\end{equation}
\noindent Remark that $\mathbf{C}_2$ is right-hand invertible, that is, there exists $\mathbf{C}_2^{-1}$ such that $\mathbf{C}_2\cdot \mathbf{C}_2^{-1}=\mathbf{Id}_{3}$. This property will be used in the control synthesis.

Only the load disk angle
position ($\theta_1$) and the main pendulum angle position ($\theta_2$)
are available by the experimental platform, but velocities will be
required in the control design. So, an observer is constructed to
estimate the velocities $\dot{\theta_1}$ and $\dot{\theta_2}$ \cite{BN93} .

\subsection{Control objective}
The control aim of this work is to design a robust control verifying
two properties. One is to ensure local stability. The other is the requirement imposed to the control design to
be robust in front of $\mathcal{L}_2$ disturbance. The
problem of robust controller with guaranteed $H_\infty$
performance is addressed to answer this question: Does there exist a feedback control such that the
$H_\infty$ norm of the closed-loop system from input disturbance
named $w(t)$ to output $z(t)$ is less than some prescribed value
$\gamma$? \cite{doyle}. In order to solve this problem, the LMI techniques are used.

 Let us consider system (\ref{ein})-(\ref{zwei}). By model variables stated previously, the vector state $\mathbf{x}^T(t)$ and the output $\mathbf{y}(t)$ are defined in (\ref{xx}) and (\ref{out}), respectively. Consider $\mathbf{z}^T(t)=[\theta_1(t)\quad
\theta_2(t)\quad u(t)]$ the virtual output to be compared to the $\mathcal{L}_2$ perturbation
$w(t)$ (induced by the second pendulum). Then, the state-space representation of
system (\ref{ein})-(\ref{out}) yields

\begin{equation}\left\{\begin{array}{l}
\dot{\mathbf{x}}(t)=\mathbf{A}\mathbf{x}(t) +
\mathbf{B}_2u(t) +
\mathbf{B}_1w(t)\;,\\[2ex]
\mathbf{z}(t)=\mathbf{C}_1\mathbf{x}(t)+\mathbf{D}_{12}u(t)\;,\\[2ex]
\mathbf{y}(t)=\mathbf{C}_2\mathbf{x}(t)\;,\\[2ex]
\end{array} \right. \label{sis3}
\end{equation}\noindent where
$$
\mathbf{A}=\left[
\begin{array}{cccc}
0&1&0&0\\[2ex]
0&-\bar{J}_z&-m^2l^2_{cg}R_hg&0\\[2ex]
0&0&0&1\\[2ex]
0&mR_hl_{cg}&ml_{cg}g(\bar{J}_1+J_y)&0
\end{array} \right] \;,
$$
$$ \mathbf{B}_1= \mathbf{B}_2=\left[
\begin{array}{c}
0\\[2ex]
\bar{J}_z\\[2ex]
0\\[2ex]
-mR_hl_{cg}
\end{array} \right],$$\vspace{0.3cm}
$$
\mathbf{C}_1=\left[
\begin{array}{cccc}
1 & 0 & 0& 0 \\
0 & 0 & 1& 0\\
0 & 0 & 0& 0 \\
\end{array} \right]  \;\textrm{and}\;
\mathbf{D}_{12}=\left[
\begin{array}{c}
0\\
0\\
1 \\
\end{array} \right].
$$

Our control objective is to design an output control matrix
$\mathbf{K}$ such that the controller
\begin{equation}\label{contr}u(t)=\mathbf{Ky}(t)\end{equation} \noindent stabilizes
the system (\ref{sis3}) under $\mathcal{L}_2$ disturbances,
employing $H_{\infty}$-LMI theory. A practical way to solve this
problem is to consider a Lyapunov function $V(\mathbf{x}(t))$ such
that for any nonzero $\mathbf{x}(t)$ and input $w(t) \in
\mathcal{L}_2$, the following condition holds:
\begin{equation}\label{dV}
\frac{d}{dt}V(\mathbf{x}(t))+\gamma^{-1}\;\mathbf{z}^T(t)\mathbf{z}(t)-\gamma
\; w^T(t)w(t) < 0\;.
\end{equation}
\noindent Then, an $H_\infty$ performance bound for the
closed-loop system (\ref{sis3})--(\ref{contr}) is ensured (see
\cite{c1} for details).
\begin{def}
\label{def-cont} If there exists a matrix $\mathbf{K}$ such that
(\ref{dV}) holds, then the control law $u(t)=\mathbf{Ky}(t)$ is
said to be an $H_{\infty}$ controller for the system
(\ref{sis3}), that is, the system is internally stable with
$H_{\infty}$ norm less than $\gamma$, i.e.,
$\|\mathbf{z}\|_{\infty}\leq \gamma^2 \|w\|_{\infty}$ for $w \in
\mathcal{L}_2$.
\end{def}

From equations (\ref{sis3}) and (\ref{contr}), we obtain the closed-loop system:
$$\left\{
\begin{array}{l}
\dot{\mathbf{x}}(t)=(\mathbf{A}+\mathbf{B}_2\mathbf{K}\mathbf{C}_2)\mathbf{x}(t) +
\mathbf{B}_1w(t)\;,\\[2ex]
\mathbf{z}(t)=\mathbf{C}_1\mathbf{x}(t)+\mathbf{D}_{12}u(t)\;.\end{array} \right.
$$
\noindent Imposing the $H_{\infty}$ condition (\ref{dV}), the term $\mathbf{B}_2\mathbf{K}\mathbf{C}_2$ difficults the evaluation of the control gain matrix. In our case, we use that $\mathbf{C}_2$ is a right-hand invertible matrix to reduce the output feedback problem with not complete information, to a full information problem \cite{doyle}:
$$
\mathbf{u}(t)=\underbrace{\mathbf{K}\mathbf{C}_2}_{\mathbf{\bar{K}}}\mathbf{x}(t) \Rightarrow \mathbf{u}(t)=\mathbf{\bar{K}}\mathbf{x}(t)\;.
$$
\noindent We solve then our problem for $\mathbf{\bar{K}}$ using the LMI theory, deriving an LMI procedure for the auxiliary gain matrix  $\mathbf{\bar{K}}$.

\emph{ Theorem 1.
Consider the Furuta pendulum system
(\ref{sis3})-(\ref{contr}). If there exist $\gamma>0$, matrices
$\mathbf{N}$, $\mathbf{Y}>0$ symmetric and $\mathbf{V}$ regular
such that the LMI
\begin{equation}\label{lmi-cond} 
\footnotesize{\left[
\begin{array}{cccccc}
-(\mathbf{V}^T+\mathbf{V})&*&*&*&*&*\\[2ex]
\mathbf{AV}+\mathbf{Y}+\mathbf{B}_2\mathbf{N}&-\mathbf{Y}&*&*&*&*\\[2ex]
0&\mathbf{B}_1^T&-\gamma&*&*&*\\[2ex]
\mathbf{C}_1\mathbf{V}&0&0&-\gamma&*&*\\[2ex]
\mathbf{N}&0&0&0&-\gamma&*\\[2ex]
\mathbf{V}&0&0&0&0&-\mathbf{Y}
\end{array}
\right]<0\;,}
\end{equation}
\noindent is feasible, then the inequality (\ref{dV}) holds, with
Lyapunov function defined as $V(\mathbf{x})=\mathbf{x}^T\;\mathbf{P}\;\mathbf{x}$ with
$\mathbf{P}:=\mathbf{Y}^{-1}$, and  $\mathbf{\bar{K}}:=\mathbf{N}\mathbf{V}^{-1}$. Consequently,
$u(t)$ is an $H_{\infty}$ controller defined by:
\begin{equation}\label{kk}
u(t)=\mathbf{N}\mathbf{V}^{-1}\mathbf{C}_2^{-1}\;\mathbf{y}(t)\;.
\end{equation}}

The proof of this result is straightforward from \cite{PA10}.

\section{Experimental realization}
To test the obtained theoretical results applied to the new laboratory device, the controller effectiveness is
studied experimentally. We design the control (\ref{contr}) via the resolution of the LMI state in {\it Theorem 1}. Experiments have been performed on an ECP Model 220 industrial
emulator with Furuta pendulum that includes a PC-based platform
and DC brushless servo system \cite{M220}.

The mechatronic system includes a motor used as servo actuator, a power
amplifier and two encoders which provide accurate position
measurements; i.e., $4000$ lines per revolution with $4X$ hardware
interpolation giving $16000 counts$ per revolution to each
encoder; $1$ count (equivalent to $0.000392$ radians  or $0.0225$
degrees) is the lowest angular measurable \cite{M220}. The second pendulum (external disturbance) includes an angular position sensor measured in radians. The drive and load disks were connected via a $4:1$ speed reduction (Figure \ref{maq-pendu}).

In the experiments, the pendulum is set to the
following parameters: $y_r=42\;cm$, $y_m=32\;cm$ (Figure \ref{441}). The parameters were taken from \cite{M220}:
$$
\mathbf{A}=\left[
\begin{array}{cccc}
0 & 1 & 0& 0 \\
0&-1.1379 & -28.769 &  0\\
0 & 0 & 0& 1 \\
0 & 0.7219& 50.229&0\\
\end{array} \right] \,,$$
$$
\mathbf{B}_1= \mathbf{B}_2=\left[
\begin{array}{c}
0 \\
318.7\\
0 \\
-202.2\\
\end{array} \right].
$$

Note that the perturbation and the control law enter to the system by the same channel ($B_{1}=B_{2}$). Using the technique presented in {\it Theorem 1}
and solving (\ref{lmi-cond}) with the LMI Matlab Toolbox \cite{matlab}, the next
$H_{\infty}$ controller was obtained: $\mathbf{\bar{K}}=\left[0.38\qquad 0.43\qquad 6.38\qquad
1.09\right]$. So, the output control (\ref{kk}) is defined as: 
\begin{equation}
\label{Kcont}
u(t)=\left[0.43\qquad 6.38\qquad1.09\right]\mathbf{y}(t)\;,
\end{equation}
\noindent with $H_{\infty}$ gain $\gamma = 8.2$. This
expression may be used directly for control modeling, scaled by appropriate system gains (amplifier and
software gains and motor torque constants). In the experiments,
the controller in (\ref{Kcont}) was multiplied by $0.3$  to
compensate these gains.

Because velocity measurements are not available, the following
controller realization is developed, where the velocity part is replaced by a first-order linear compensator:
\begin{equation}\label{cont-alg}
\begin{array}{lll}
u=0.3\left(0.43\;\dot{x}_1+6.38\;\theta_2+1.09\;\dot{x}_2\right)\;,\\
\dot{x}_1=-10x_1+5\theta_1\;,\\\dot{x}_2=-10x_2+5\theta_2\;.
\end{array}
\end{equation}

The terms $x_1$ and $x_2$ are auxiliary variables and their differential equations are solved numerically from $\theta_1$ and $\theta_2$ measurements.  The above equations are a direct implementation of
velocity observers given in \cite{BN93}, where the parameter $-10$
was set according to \cite{Ogata} and the gain $5$
was adjusted experimentally\footnote{This consist to locate
the observer poles from $5$ to $10$ times far away with respect to
the vertical-axis-closest pole of $(A+B_2K)$.}. This approach goes along a two independent steps
design procedure: a) design an output state-feedback controller $\mathbf{K}$, and b) construct a
velocity observer. This design
obeys the separation principle (see \cite{Ogata} and \cite{K00}).\\

In this experiment we try to simulate the Segway navigation, studying the load disk displacement when the main pendulum is strongly perturbed. Figure \ref{fe1a} presents a traveling from the video http://youtu.be/SHQdW3k3qiE, where the complete navigation experiment can be seen. Figures \ref{fe3}-\ref{fe5} show  the load and principal pendulum displacement (as Segway navigation), and the control effort. Remark that at $t=8$ sec, we return the second pendulum to its  initial position (as Figure \ref{fe1a} shows). Then, the load disk remains in its final position (the Segway does not move if the person comes back to its upright position). This phenomena can be seen in Figure \ref{fe3}.
\begin{figure}[ht!]
\begin{center}
\caption{Clockwise direction. Pictures from http://youtu.be/SHQdW3k3qiE.}  
\label{fe1a}                                 
\end{center}                                 
\end{figure}
\begin{figure}[ht!]
\begin{center}
\caption{Load disk position (simulation of the Segway navigation).}  
\label{fe3}                                 
\end{center}                                 
\end{figure}
\begin{figure}[ht!]
\begin{center}
\caption{Main pendulum position (simulation of driving the Segway).}  
\label{fe4}                                 
\end{center}                                 
\end{figure}
\begin{figure}[ht!]
\begin{center}
\caption{Control effort.}  
\label{fe5}                                 
\end{center}                                 
\end{figure}

\section{Conclusion}
A new experimental set up for the Furuta pendulum was developed to validate controller performance under a kind of mass perturbation similar to the one used to navigate the Segway personal transportation unit. To simulate the Segway behavior, a second pendulum was coupled to the main one. The output control was designed in order to ensure stability and robustness. Experimental results demonstrate that the objectives have been achieved. Hence, the mechanical experiment provides interesting new results in the control automatic new developments. 

\section{Acknowledgment}
This work was partially funded by the Spanish Government, the European Commission, and  FEDER funds, through the research projects: CGL2008-00869/BTE, CGL2011-23621, DPI2011-26326, DPI2011-25822 and DPI2012-32375/FEDER.
%
%

%
%
%
%

\end{document}